# Proof of Goldbach's Conjecture

Reza Javaherdashti
farzinjavaherdashti@gmail.com

**Abstract**

After certain subsets of Natural numbers called "Range" and "Row" are defined, we assume (1) there is a function that can produce prime numbers and (2) each even number greater than 2, like A, can be represented as the sum of n prime numbers. We show this by $DC(A) \leq n$. Each Row is similar to each other in properties,(so is each Range). It is proven that in an arbitrary Row for any even number greater than 2, $DC(A)=2$, that is to say, each even number greater than two is the sum of two prime numbers. So Goldbach's conjecture is proved.

## 1.Historical Background:

Of still-unsolved problems on prime numbers one can mention Goldbach's conjecture. Goldbach (1690-1764) in his letter to Euler in 1742, asked if any even number greater than 2 could be written as the sum of two prime numbers. Euler could not answer nor could he find any counter-example. The main problem with Goldbach's conjecture is that in most of theorems in arithmetic, prime numbers appear as products, however, in Goldbach's conjecture it is the addition of prime numbers that makes all the problem. In 1931, a young, not that famous Russian mathematician named Schnirelmann (1905-1938) proved that any positive integer could have been represented, at most, as the sum of 300,000 prime numbers. The reasoning was constructive and direct without giving any practical use to decompose a given integer into the sum of prime numbers. Some years after him, the Russian mathematician Vinogradoff by using and improving methods invented by the English mathematicians, Hardy and Littlewood, and their great Indian colleague Ramanujan, could decrease number of the mentioned prime numbers from 300,000 to 4. Vinogradoff's approach has been proved to be true for integers "large enough". With exact words, Vinogradoff proves that there exists an integer like N so that for any integer n>N, it can be represented, at most, as the sum of four prime numbers. He gives noway to determine and measure N. Vingoradoff's method has actually proved that accepting the infinite integers that cannot be shown as the sum of, at most, four prime numbers, results in contradiction.

**2.Method and Basic Assumptions:**

The method we use in our approach, relies on the following facts (or, interpretations):
1.Controversial points with Goldbach's conjecture are:

> 1.1.It seems as if there must be a kind of formula that can produce prime numbers.
> 1.2. After such a prime number-producing formula is found, one should look for its relationship with even numbers greater than 2.

2.Goldbach assumes that sum of prime numbers gives even numbers greater than 2; the problem is how to limit the number of such prime numbers with only two.
To clarify the approach, we assume that:
*Assumption #1*: There is a function like f (x) that produces prime numbers.
*Assumption #2*: Each even number greater than 2, can be taken as to be the sum of n prime number where n is a Natural number.

Using the above-mentioned assumptions, we define "Row" and "Range" as subsets of Natural numbers. Then, by using DC (A), that designates that how many prime numbers can produce even number "A" greater than 2, we prove DC(A)=2. In other words, It is proved that that the minimum number of prime numbers that results in an even number "A", which is greater than 2, is only two. Thus, Goldbach's conjecture is proved.
Our method consists of three parts:
I.    Basic concepts on "Row" and "Range".
II.   Basic definitions of f(x) and DC(A).
III.  Proof.
We will NOT use statistical data or tables of prime numbers in our method.

## Part 1): Basic concepts on "Row" and "Range"

### 1.1. Definition of Row:

Row-that we show as $r(x_i\ x_f)$- is a term used for representing any subset of Natural numbers, N, that has <u>all</u> of the properties below:

I.  $r(x_i\ x_f) \subset N$, that is, each "Row" is a subset of Natural numbers.
    If $n[r(x_i\ x_f)]$ shows the number of elements of the set $r(x_i\ x_f)$, then :
    $$n[r(x_i\ x_f)] = d \qquad d \in N$$
    The above means that in a Row, number of elements is limited.

II. (property of having a smallest and a greatest element in a Row):
    Each $r(x_i\ x_f)$ has just one "smallest" and just one "greatest" element, that is, each $r(x_i\ x_f)$, at most, has one "smallest" element like $x_i$ in $r(x_i\ x_f)$ such that for all $x \in r(x_i.x_f)$, $x > x_i$. In the same way, each $r(x_i\ x_f)$ has one "greatest" element like $x_f$ such that for all $x \in r(x_i\ x_f)$, $x < x_f$. So $x_f$ is the "greatest element" of $r(x_i\ x_f)$ and there is no element larger than it.

III. (property of an ordered Row):
    In each Row $r(x_i\ x_f)$ all of its elements are orderable, from left to right, and from the smallest to the largest element.

IV. (property of constant difference between elements of a Row)
    For all $x \in r(x_i\ x_f)$, if $x \neq x_i$ and $x \neq x_f$, then $x-1 < x < x+1$. This means that each $x \in r(x_i\ x_f)$ is less than or larger than a number immediately before or after by a constant difference, which is one.

*1.1.1. Examples:*

A = {1,2,3,4}   A is a Row  so A = r(14) as $x_i = 1$ and $x_f = 4$ and the difference between each immediate successive element is unity.

B = {25,26,27} B is a Row so B = r(2527) as $x_i = 25$ and $x_f = 27$

C = {4,3,2,1}    C is NOT a Row; III is not held.

D = {5,9,10,11,14}   D is NOT a Row; IV is not held.

E = {49,51,53,55}    E is NOT a Row; IV is not held.

**Convention.1:**

1. From now on, instead of $r(x_i \, x_f)$, the symbol $r(if)$ is used.
2. According to III, $x < x_i$ or $x > x_f$ is not defined in $r(if)$.
3. Any $r(if)$, schematically can be shown as follows:

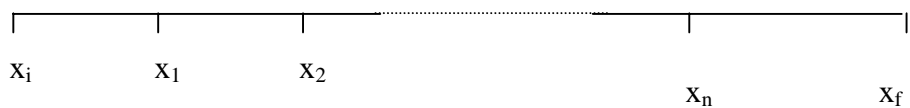

$$r(if) = \{ x_i, x_1, x_2, \ldots, x_n, x_f \} \qquad x_n < x_f$$

## 1.2. Definition of Range:

Range-that we show as $R(x_I \, x_F)$-is a term used for representing any subset of Natural numbers that has <u>all</u> of the properties below:

[1]. $R(x_I \, x_F) \subset N$

[2]. $r(if) \subset R(x_I \, x_F)$

[3]. $n[R(x_I \, x_F)] = D \qquad D \in N, \quad D > d$

[4]. Each $R(x_I \, x_F)$ has, at most, one "smallest" element shown as $x_i$, and one "greatest" element $x_F$, ie; for each $X \in R(x_I \, x_F)$:

$$X < x_F, \, X > x_I$$

[5]. For each $X \in R(x_I \, x_F)$ giving $X \neq x_I$ and $X \neq x_F$:

$$X-1 < X < X+1$$

[6]. All the elements of a Range $R(x_I \, x_F)$ can be ordered, from left to right and from the smallest to the largest element.

*1.1.2. Examples:*

$F = \{1,2,3,\ldots,98,99,100\}$    F is a Range so $F = R(1100), r(110) \subset F,\ldots$ so [2] is held.

$G = \{5,6,7,\ldots,22,23\}$    G is range so $G = R(523)$

$H = \{30,31,32,35,36,37,39,41,42,45\}$    H is NOT a Range, [5] is not held.

$I = \{31,32,50,33,50,1,2,4,10,2000\}$

   I, is NOT a Range,;

   -[2] is not held as only $\{1,2\}$ and $\{31,32\}$ have the smallest and the largest elements, so IV is not held for these subsets.

-[5] is not held.

-[6] is not held.

**Convention.2:**

1-From now on, $R(x_I\ x_F)$ is shown as R(IF).

2-By considering [2] and [4], it appears that in some cases $X_I = x_i$ or $X_F = x_f$ ,for example, in F, $X_I = 1$ and also in A, $x_i = 1$ (A $\subset$ F).

3-By using "Row" and "Range" concepts, N (natural numbers set) can be subdivided into subsets that have certain properties as R(IF) $\subset$ N and r(if) $\subset$ R(IF). There are infinite number of R(IF) sets but definite and limited number of r(IF) subsets, for instance, there are inifinite number of sets like R(1100), R(101200), R(201300) etc. but in each of these sets there are limited and definite number of r(if);as an example, in R(1100),there are ten subsets (=Rows) such as r(110), r(1120), r(2130),…, r(91100). So, any theorem that is proved or any conclusion that is made for arbitrary r(if) and R(IF), can then be generalised for all similar Rows and Ranges so that what is proven, will be applicable to whole Natural numbers.

**Convention.3:**

For each x which is an element of a Row r(if), we use the following symbolism:

$$x \in r(if) \text{ or } (x)r(if)$$

and for each X which is an element of R(IF), we use the following symbolism:

$$X \in R(IF) \text{ or } (X)R(IF)$$

And so forth.

**Convention.4**:

In a Row r(ij), its smallest element-that is $x_i$ –can be shown as either $(x_i)r(ij)$ or A( i ) and its greatest element-$x_f$-can be shown as $(x_f)r(ij)$ and so forth.

**Result of Convention.4:**

In a Row r(pq), A(p) is the smallest element. In a Row r(st), A(s) is the smallest element and so forth.

*Theorem.1*

Take two Rows r(ij) and r(kl) in a range R(IF). Prove that for all $x_i$ and $x_f$ elements of these Rows, we have:

$$(x_f)r(ij) \pm p = (x_i)r(kl)$$

where p is a constant.

*Proof:*

Assume that p = 0, then $(x_f)r(ij) = (x_i)r(kl)$. Now take R(IF) as:

(1) $(r(ij) \subset R(IF))$ & $(r(kl) \subset R(IF)) \leftrightarrow R(IF) = r(ij) \cup r(kl)$, where "∪" is union sign, results in:

(2) $(x \in R(IF)) \leftrightarrow (x \in r(ij))$ & $(x \in r(kl))$

as $(x_i)r(kl) \in r(kl)$ and $(x_f)r(ij) \in r(ij)$, then both $(x_i)r(kl)$ and $(x_f)r(ij)$ must be elements of R(IF) but in this case:

1. the difference between $(x_i)r(kl)$ and $(x_f)r(ij)$ will not be equal to unity ( [5] and IV are not held).

2. R(IF) will not be ordered ( [6] is not held)

So, R(IF) will not be a Range but this is in contradiction with our assumption. Therefore, p≠0. If each element of r(ij) is smaller than each element of r(kl) then, $(x_f)r(ij) + p = (x_i)r(kl)$ and so forth.

**Definition.1**:

If $(x_f)r(ij) + p = (x_i)r(kl)$, r(ij) and r(kl) are called "Successive Rows".

**Reminder:**

"Range"s are arbitrary, random sample sets from set of natural numbers. "Row"s are arbitrary, random sample sets from ranges. By introducing the concepts of range and row, it is followed that if something is proven for these random "cuts" from natural numbers set, it has been proven for WHOLE natural numbers. However, definitions of range and row clearly show that NOT ANY random set from natural numbers set can be picked up as range and/or row (see examples 1.1.1. and 1.1.2.).

**Part.2): Basic definitions of f(x) and DC(A)**

**Definition 2:**

We define a function f(x) so that for any x = a, f(a) be a prime number.

**Definition 3:**

The "degree of complexity" function that we show it as "DC" is the number of prime numbers to be added to each other with a + sign between each to yield an even number greater than two. DC will itself be a Natural number.

*2.2.1.Example*:

-For 8 = 2 + 2 + 2 + 2         DC (8) = 4

-For 8 = 5 + 3                    DC (8) = 2

-For 216 = 213 + 3                DC (216) = 2

**Conclusion from definition 3:** The least value for DC(A), where A is an even number greater than 2, is 2; ie, DC (A) ≥ 2.

**Definition 4:**

In any Row r(ij) of the Range R(IF), number of even numbers (2K), odd numbers (2K + 1) and prime numbers (f(x)) are shown, respectively, as Γ(2K), Γ(2K + 1), and Γf(x). So in the Range R(1100) and Row r(110), Γ(2K) = 5 that is to say, (2,4,6,8,10), Γ(2K + 1) = 5 ie (1,3,5,7,9) and Γf(x) = 4 (2,3,5,7).

**Part.3): Proof**

Assumptions:

1) Assume in Row r(ij) of Range R(IF), there exists at least one prime number like f(x).

2) Assume in any Row r(ij), the number of even and odd numbers are equal to each other (in any Row there is as many odd numbers as there is even numbers). So in any Row r(ij) there exists odd and even numbers alternatively (after each odd number there is an even number and vice versa).

3) Assume in a Row r(ij), number of odd numbers be more than prime numbers in the same Row, in other words, any prime number greater than 2 is an odd number BUT any odd number is NOT a prime number.

From assumptions 1) and 3), one concludes:

(3) $1 \leq \Gamma f(x) \leq \Gamma(2K + 1)$

For even numbers greater than 2, like A, degree of complexity function can be written as DC (A). According to the conclusion from definition 3, DC (A) ≥ 2, ie, at least two prime numbers must be added to each other to yield A. As A is an even number greater than 2, the number of even numbers that are required to be added to each other to yield A, will be less than the number of existing even number in the Row r(ij), so in r(ij) the number of prime numbers to be added up to yield A is:

(4) $DC (A) \leq \Gamma(2K)$

By adding Γ(2K + 1) to right-hand sides of inequalities (3) and (4) it yields:

(5) $\Gamma f(x) \leq \Gamma(2K + 1) + \Gamma(2K + 1)$

(6) DC (A) ≤ $\Gamma(2K) + \Gamma(2K + 1)$

According to assumption 2) $\Gamma(2K) = \Gamma(2K + 1)$; by applying this to (5), it yields:

(7) $\Gamma f(x) \leq \Gamma(2K) + \Gamma(2K + 1)$

by adding each side of (6) and (7) to each other, we take:

(8) DC (A) + $\Gamma f(x) \leq 2[\Gamma(2K) + \Gamma(2K + 1)]$

or

(9) DC (A) ≤ $2[\Gamma(2K) + \Gamma(2K + 1)] - \Gamma f(x)$

Combining the conclusion from definition 3 with inequality (9) yields:

(10)    $2 \leq$ DC (A) $\leq 2[\Gamma(2K) + \Gamma(2K + 1)] - \Gamma f(x)$

The above relation can be decomposed into the following three inequalities:

(11-1)  DC (A) > 2

(11-2)  DC (A) $\leq 2[\Gamma(2K) + \Gamma(2K + 1)] - \Gamma f(x)$

(11-3)  DC (A) = 2

We will prove that (11-1) and (11-2) will be resulting in contradictions so that they will not be held. Therefore, the only remaining relation will be (11-3) that states that the number of prime numbers to be added up to yield an even number greater than 2, is two.

Assume (11-1) holds, ie, DC > 2. This means:

(12)   DC (A) = $a_1 + a_2 + \sum_{i=3}^{i=m} a_i$

Where m shows total number of numbers-even, odd and prime numbers-existing in a given Row r(ij); $a_i$ (i = 1,2,3,…,m) shows the prime numbers to yield A. Relation (12) may be re-written as (13):

(13)   DC (A) = $2 + \sum_{i=n}^{i=m} a_i$

Where n is an arbitrary number less than m. One should notice that in (12) $a_1$ and $a_2$ are two even numbers where in (13) some prime numbers like $\sum_{i=n}^{i=m} a_i$ are added to a number like 2. Equation (13) yields:

(14)   DC (A) – 2 = $\sum_{n}^{m} a_i$

$\sum_{n}^{m} a_i$ represents the number of prime numbers in a Row r(ij) which is less than total number of existing prime numbers of r(ij), that is to say:

(15) $\sum_{n}^{m} a_i \leq \Gamma f(x)$

replacing (14) into (15) yields:

(16) DC (A) – 2 = $\Gamma f(x)$

Assumption 3) of **Part.3):Proof** about odd and prime numbers gives:

(17) $\sum_{n}^{m} a_i \leq \Gamma(2K + 1)$

replacing (14) into (17) yields:

(18) DC (A) – 2 = $\Gamma(2K + 1)$

Or

(19) DC (A) $\leq \Gamma(2K + 1) + 2$

From (16), (20) is resulted:

(20) DC (A) $\leq \Gamma f(x) + 2$

By adding sides of (19) and (20) to each other:

(21) 2DC (A) $\leq \Gamma(2K) + \Gamma f(x) + 4$

as DC (A) < 2DC (A), one may conclude (22). To let the inequality hold with more force, we add $\Gamma(2K)$ to right-hand side of (21) too;

(22) 2DC (A) $\leq \Gamma(2K) + \Gamma(2K + 1) \Gamma f(x) + 4$

The right-hand side of (22) may be written as (23):

(23) $2[\Gamma(2K) + \Gamma(2K + 1)] - \Gamma f(x) = [\Gamma(2K) + \Gamma(2K + 1) + \Gamma f(x) + 4] + [\Gamma(2K) + \Gamma(2K + 1) - 2\Gamma f(x) - 4]$

(24) $\Gamma(2K) + \Gamma(2K + 1) + \Gamma f(x) + 4 < 2[\Gamma(2K) + \Gamma(2K + 1)] - \Gamma f(x)$

To combine (22) and (24) yields:

(25) DC (A) $< \Gamma(2K) + \Gamma(2K + 1) + \Gamma f(x) + 4 < 2[\Gamma(2K) + \Gamma(2K + 1)] - \Gamma f(x)$

The relation (25) results in:

(26) DC (A) $< 2[\Gamma(2K) + \Gamma(2K + 1)] - \Gamma f(x)$

Inequality (26) resembles (11-2). So we will consider (26) more precisely:

One sees that DC (A) $\leq 2[\Gamma(2K) + \Gamma(2K + 1)] - \Gamma f(x)$ shows the number of prime numbers in a given Row r(ij) to be added to each other to result in a number like A, which is an even number greater than 2. Therefore, $2[\Gamma(2K) + \Gamma(2K + 1)] - \Gamma f(x)$ must be less than the number of existing odd numbers in r(ij) (Assumption 3) of **Part.3):Proof).**

On the other hand, according to assumption 2) **Part.3):Proof** at least half of the total numbers in a Row r(ij) –in other words, $\frac{m}{2}$ -are existing odd numbers of r(ij). As it was stated in (12), if m be total number of existing numbers in a Row r(ij), then:

(27) $2[\Gamma(2K) + \Gamma(2K + 1)] - \Gamma f(x) \leq \frac{m}{2}$

as $\Gamma(2K) + \Gamma(2K + 1) = m$ (Assumption 2) of **Part.3):Proof)):**

(28) $2m - \Gamma f(x) \leq \frac{m}{2}$

(29) $1.5m \leq \Gamma f(x)$

as m < 1.5m then:

(30) $m < 1.5m < \Gamma f(x)$

or:

(31) $m < \Gamma f(x)$

Inequality (31) is not held as it states that in a Row r(ij), the number of prime numbers of the Row is larger than total numbers of the Row, which is IMPOSSIBLE. In the same way, one may prove that if (26) is written as:

DC (A) = $2[\Gamma(2K) + \Gamma(2K + 1)] - \Gamma f(x)$, again we will come up with the contradiction stated as in (31).

The above discussion shows that (11-2) is NOT held.

One may show (25) as a combination of the following:

(32-1) DC (A) > 2

(32-2) DC (A) < $\Gamma(2K) + \Gamma(2K + 1) \Gamma f(x) + 4$

(32-3) DC (A) < $2[\Gamma(2K) + \Gamma(2K + 1)] - \Gamma f(x)$

We call (32-1) through (32-3), as p, q and r, respectively so that we translate (25) to the logical expression (33)as a conjunction:

(33) p $\implies$ ( q & r)

where "&" represents conjunction sign. According to what we have gained so far::

1.r is false, ie, relations (27) to (31) prove that :

DC (A) $\leq 2[\Gamma(2K) + \Gamma(2K + 1)] - \Gamma f(x)$, is NOT held.

2.as r is false, then (q & r) is false. So p must be false also , that is to say:

DC (A) > 2 is NOT held.

**Conclusion**
01. (11-1) isnot held so DC (A) > 2 isnot true.
02. (11-2) isnot held so DC (A) ≤ 2[Γ(2K) + Γ(2K + 1)] - Γf(x) isnot true.
03. (11-3) holds as the only possibility so DC (A) = 2 is true.

In fact Conclusion 03 states that an even number greater than 2 can be written as the sum of two prime numbers. So Goldbach's conjecture is true.

**Summary:**

By defining certain subsets in Natural numbers set, after all possibilities are considered, we prove that Goldbach's conjecture is held in one arbitrary. The proof can be generalised to all the subsets so that the conjecture is proven for all Natural numbers.